\documentclass[a4paper,3p]{elsarticle}
\usepackage{mathtools}
\usepackage[pdftex,colorlinks,linkcolor=blue,citecolor=blue,urlcolor=blue]{hyperref}
\usepackage{amssymb}
\usepackage{amsthm}
\usepackage{stix2}
\usepackage{inconsolata}
\usepackage{bm}
\usepackage{mleftright}
\mleftright
\usepackage{algorithm}
\usepackage{algorithmic}
\numberwithin{equation}{section}

\def\rme{\mathrm{e}}
\def\rmi{\mathrm{i}}
\def\Rii{$\mathrm{R}_{\mathrm{II}}$\ }

\theoremstyle{definition}
\newtheorem*{example}{Example}

\newcommand{\IfThen}[2]{
  \algorithmicif\ #1\ \algorithmicthen\ #2}
\newcommand{\IfThenElse}[3]{
  \algorithmicif\ #1\ \algorithmicthen\ #2\ \algorithmicelse\ #3}

\begin{document}

\title{An isospectral transformation between Hessenberg matrix and Hessenberg--bidiagonal matrix pencil without using subtraction}
\author[1]{Katsuki Kobayashi}
\author[2]{Kazuki Maeda\corref{cor1}}
\ead{kmaeda@kmaeda.net}
\author[1]{Satoshi Tsujimoto}
\ead{tsujimoto.satoshi.5s@kyoto-u.jp}
\cortext[cor1]{Corresponding author}
\affiliation[1]{organization={Department of Applied Mathematics and Physics, Graduate School of Informatics, Kyoto University},
  city={Yoshida-Honmachi},
  postcode={Kyoto 606-8501},
  country={Japan}}
\affiliation[2]{organization={Faculty of Informatics, University of Fukuchiyama},
  city={3370 Hori, Fukuchiyama},
  postcode={Kyoto 620-0886},
  country={Japan}}

\begin{abstract}
We introduce an eigenvalue-preserving transformation algorithm from
the generalized eigenvalue problem by matrix pencil of the upper and
the lower bidiagonal matrices into a standard eigenvalue problem while preserving sparsity,
using the theory of orthogonal polynomials. The procedure is formulated without subtraction,
which causes numerical instability. Furthermore, the algorithm is discussed
for the extended case where the upper bidiagonal matrix is of Hessenberg type.

\end{abstract}


\maketitle

\section{Introduction}

Generalized eigenvalue problems (GEVP) appear in various fields and are often attributed to
standard eigenvalue problems (EVP). If the matrix $B$ of the matrix pencil $(A,B)$ is regular,
then it is possible to consider its inverse and solve it as an EVP for the matrix $B^{-1}A$.
However, even when both $A$ and $B$ are sparse matrices, the matrix $B^{-1}A$ will lose
the sparsity of $A$ and $B$. This will lead to computational complexity and numerical instability.
Hence it is useful to consider a transformation from GEVP to EVP preserving the sparsity of
the original matrix without using subtractions, which can cause numerical instability,
if possible.

The simplest GEVP treated here is a matrix pencil consisting of a tridiagonal matrix
$A$ and an almost diagonal matrix $B$, in which a nonzero component appears in the subdiagonal
in one place:
\begin{gather*}
  A=
  \begin{pmatrix}
    a_{0, 0} & a_{0, 1}\\
    a_{1, 0} & \ddots & \ddots\\
    & \ddots & a_{j, j} & a_{j, j+1}\\
    && a_{j+1, j} & a_{j+1, j+1} & \ddots\\
    &&& \ddots & \ddots & a_{N-2, N-1}\\
    &&&& a_{N-1, N-2} & a_{N-1, N-1}
  \end{pmatrix},\quad
  B=
  \begin{pmatrix}
    \alpha_0\\
    & \ddots\\
    && \alpha_j\\
    && \beta & \alpha_{j+1}\\
    &&&& \ddots \\
    &&&&& \alpha_{N-1}
  \end{pmatrix},
\end{gather*}
where $\beta\ne 0$.
In this case, we first convert the GEVP $A\bm\phi = \lambda B\bm\phi$ to the following GEVP
\begin{align*}
 \left(A-\frac{a_{j+1,j}}{\beta}B\right)\bm\phi = \lambda' B\bm\phi,
\end{align*}
where $\lambda'=\lambda-\frac{a_{j+1,j}}{\beta}$.
We note here that in the tridiagonal matrix $A-\frac{a_{j+1,j}}{\beta}B$,
one of the subdiagonal components is a zero component.
We will show that GEVP in the above form can be rewritten as an EVP without subtraction
by using the theory of Laurent biorthogonal polynomials~\cite{kharchev1997frt,zhedanov1998clb}.
Furthermore, the same procedure can be introduced for the generalization of the above form.
The following is an outline of the concrete forms and procedures of the GEVP treated in this paper.

In Section~\ref{sec:transf-betw-trid},
first, we consider a simple GEVP
\begin{equation*}
 R \bm p = x L \bm p
\end{equation*}
of the pencil $(R, L)$ defined by the upper bidiagonal matrix $R$
and the lower bidiagonal matrix $L$ as
\begin{equation*}
  R\coloneq
  \begin{pmatrix}
    q_{0} & 1\\
    & q_{1} & 1\\
    && \ddots & \ddots\\
    &&& \ddots & 1\\
    &&&& q_{N-1}
  \end{pmatrix}, \quad
   L\coloneq
  \begin{pmatrix}
    1 \\
    -e_0 & 1\\
    & -e_1 & \ddots\\
    && \ddots & \ddots\\
    &&& -e_{N-2} & 1
 \end{pmatrix}.
\end{equation*}
We can see that each element of the vector $\bm p$ is a polynomial of
$x$ defined by a three-term recurrence relation for
the Laurent biorthogonal polynomials.
The followings will be shown:
\begin{itemize}
\item LU factorization of the matrix $RL^{-1}$
  yields new bidiagonal matrices $R^*$, $L^*$ and the associated vector (polynomials) $\bm p^*$,
  where the new GEVP $R^*\bm p^*=x L^*\bm p^*$ has the same eigenvalues as
  the original GEVP.
\item Setting $(R^{(0)}, L^{(0)})\coloneq (R, L)$, $\bm p^{(0)}\coloneq \bm p$ and iterating the procedure above,
  we obtain the matrix pencil sequence $(R^{(0)}, L^{(0)})$, $(R^{(1)}, L^{(1)})$,
  $(R^{(2)}, L^{(2)})$, ...,
  and the corresponding vector sequence $\bm p^{(0)}$, $\bm p^{(1)}$, $\bm p^{(2)}$, ....
  Then, we can construct a tridiagonal matrix $\hat T$
  from the matrix pencil sequence and another vector $\hat{\bm p}$ from the vector sequence,
  where $\hat{\bm p}$ is the eigenvector of the EVP
  $\hat T\hat{\bm p}=x\hat{\bm p}$
  which has the same eigenvalues as the original GEVP.
  Elements of the vector $\hat{\bm p}$ are orthogonal polynomials~\cite{chihara1978iop} defined by
  the three-term recurrence relation whose coefficients are given by $\hat T$.
  Note that, since $L^{-1}R$ is not a tridiagonal matrix,
  this is not a trivial isospectral transformation.
\end{itemize}

In Section~\ref{sec:transf-betw-trid-1},
the transformation in Section~\ref{sec:transf-betw-trid} is generalized to
a matrix pencil $(T, L)$, where $T$ is the tridiagonal matrix
and $L$ is the lower bidiagonal matrix of the form
\begin{gather*}
  T=
  \begin{pmatrix}
    a_0 & 1\\
    (1-\epsilon_0) b_0 & a_1 & 1\\
    & (1-\epsilon_1) b_1 & \ddots & \ddots\\
    && \ddots & \ddots & 1\\
    &&& (1-\epsilon_{N-2}) b_{N-2} & a_{N-1}
  \end{pmatrix},\\
  L=
  \begin{pmatrix}
    1\\
    -\epsilon_0 e_0 & 1\\
    & -\epsilon_1 e_1 & \ddots\\
    && \ddots & \ddots\\
    &&& -\epsilon_{N-2} e_{N-2} & 1
  \end{pmatrix},\quad
  \epsilon_0, \dots, \epsilon_{N-2} \in \{0, 1\}.
\end{gather*}
It will be shown that
there is also a transformation to a tridiagonal matrix with the same eigenvalues.

In Section~\ref{sec:transf-betw-hess}, further,
the transformation in Section~\ref{sec:transf-betw-trid-1}
is generalized to a matrix pencil $(H, L)$, where
$H$ is the upper Hessenberg matrix and $L$ is the lower bidiagonal matrix
of the form
\begin{gather*}
  H=
  \begin{pmatrix}
    a_{0, 0} & a_{0, 1} & \dots & a_{0, M-1} & 1\\
    (1-\epsilon_0) b_0 & a_{1, 0} & a_{1, 1} & \dots & a_{1, M-1} & 1\\
    & (1-\epsilon_1) b_1 & \ddots & \ddots & \dots & \ddots & \ddots\\
    && \ddots & \ddots & \ddots & \dots & \ddots \\
    &&& \ddots & \ddots & \ddots & \dots\\
    &&&& \ddots & \ddots & a_{N-2, 1}\\
    &&&&& (1-\epsilon_{N-2}) b_{N-2} & a_{N-1, 0}
  \end{pmatrix},\\
  L=
  \begin{pmatrix}
    1\\
    -\epsilon_0 e_0 & 1\\
    & -\epsilon_1 e_1 & \ddots\\
    && \ddots & \ddots\\
    &&& -\epsilon_{N-2} e_{N-2} & 1
  \end{pmatrix},\quad
  \epsilon_0, \dots, \epsilon_{N-2} \in \{0, 1\}.
\end{gather*}
It will be shown that
there is a transformation to an upper Hessenberg matrix with the same eigenvalues.

\section{Isospectral transformation between tridiagonal matrix and bidiagonal--bidiagonal matrix pencil}
\label{sec:transf-betw-trid}

\subsection{Sequence of bidiagonal--bidiagonal GEVPs and eigenvectors}
Let $N$ be a positive integer and $\{p_n(x)\}_{n=0}^{N}$ be monic polynomials of degree $n$
defined by the three-term recurrence relation
\begin{gather}
  p_{-1}(x)\coloneq 0,\quad p_0(x)\coloneq 1,\quad
  p_{n+1}(x) \coloneq (x-q_n)p_n(x)-x e_{n-1} p_{n-1}(x),\quad n=0, 1, \dots, N-1,\label{eq:trr-LBP}
\end{gather}
where $q_n, e_n\in \mathbb C$ and $e_n\ne 0$.
We can rewrite the three-term recurrence relation~\eqref{eq:trr-LBP} as
\begin{equation}\label{eq:trr-LBP-vec}
  R\bm p(x)+\bm p_N(x)=xL\bm p(x),
\end{equation}
where $R$ and $L$ are bidiagonal matrices
\begin{equation*}
  R\coloneq
  \begin{pmatrix}
    q_0 & 1\\
    & q_1 & 1\\
    && \ddots & \ddots\\
    &&& \ddots & 1\\
    &&&& q_{N-1}
  \end{pmatrix},\quad
  L\coloneq
  \begin{pmatrix}
    1 \\
    -e_0 & 1\\
    & -e_1 & \ddots\\
    && \ddots & \ddots\\
    &&& -e_{N-2} & 1
  \end{pmatrix},
\end{equation*}
and $\bm p(x)$, $\bm p_N(x)$ are vectors of polynomials
\begin{equation}\label{eq:def-vec-p}
  \bm p(x)\coloneq
  \begin{pmatrix}
    p_0(x)\\
    p_1(x)\\
    \vdots\\
    p_{N-1}(x)
  \end{pmatrix},\quad
  \bm p_N(x)\coloneq
  \begin{pmatrix}
    0\\
    \vdots\\
    0\\
    p_N(x)
  \end{pmatrix}.
\end{equation}
Let $x_0, x_1, \dots, x_{N-1}$ be the zeros of $p_N(x)$, then we have
\begin{equation*}
  R\bm p(x_i)=x_i L\bm p(x_i),\quad i=0, 1, \dots, N-1,
\end{equation*}
i.e., $x_i$ and $\bm p(x_i)$ are a generalized eigenvalue and eigenvector
of the matrix pencil $(R, L)$.

Next, let us introduce new monic polynomials $\{p_n^*(x)\}_{n=0}^N$ of degree $n$ generated from
the monic polynomials $\{p_n(x)\}_{n=0}^N$ as
\begin{gather}
  p_0^*(x)\coloneq 1,\quad p_N^*(x)\coloneq p_N(x),\quad
  p_n^*(x)\coloneq p_n(x)-e_{n-1}p_{n-1}(x),\quad n=1, 2, \dots, N-1.\label{eq:GT-LBP}
\end{gather}
Then the three-term recurrence relation of $\{p_n(x)\}_{n=0}^N$ \eqref{eq:trr-LBP} is
rewritten as
\begin{equation}
  xp_n^*(x)=p_{n+1}(x)+q_n p_n(x),\quad n=0, 1, \dots, N-1.\label{eq:CT-LBP}
\end{equation}

Using the bidiagonal matrices and vectors of polynomials,
the relations \eqref{eq:GT-LBP} and \eqref{eq:CT-LBP} are rewritten as
\begin{equation}\label{eq:GT-CT-LBP}
  \bm p^*(x)=L\bm p(x),\quad
  x\bm p^*(x)=R\bm p(x)+\bm p_N(x),
\end{equation}
where
\begin{equation*}
  \bm p^*(x)=
  \begin{pmatrix}
    p_0^*(x)\\
    p_1^*(x)\\
    \vdots\\
    p_{N-1}^*(x)
  \end{pmatrix}.
\end{equation*}
Since $L$ is a regular matrix, we obtain
\begin{equation}\label{eq:next-LBP}
  x\bm p^*(x)=RL^{-1}\bm p^*(x)+\bm p_N(x).
\end{equation}
The LU factorization of $RL^{-1}$ generates bidiagonal matrices
\begin{equation*}
  R^*\coloneq
  \begin{pmatrix}
    q_0^* & 1\\
    & q_1^* & 1\\
    && \ddots & \ddots\\
    &&& \ddots & 1\\
    &&&& q_{N-1}^*
  \end{pmatrix},\quad
  L^*\coloneq
  \begin{pmatrix}
    1 \\
    -e_0^* & 1\\
    & -e_1^* & \ddots\\
    && \ddots & \ddots\\
    &&& -e_{N-2}^* & 1
  \end{pmatrix},
\end{equation*}
satisfying
\begin{equation*}
  RL^{-1}=(L^*)^{-1}R^*.
\end{equation*}
Then, substituting this relation to \eqref{eq:next-LBP},
we obtain
\begin{equation*}
  R^*\bm p^*(x)+\bm p_N(x)=x L^* \bm p^*(x).
\end{equation*}
This means that the monic polynomials $\{p_n^*(x)\}_{n=0}^N$ also satisfy
three-term recurrence relation of the same form as $\{p_n(x)\}_{n=0}^N$.
Therefore $x_i$ and $\bm p^*(x_i)$ are a generalized eigenvalue and eigenvector
of the matrix pencil $(R^*, L^*)$.
Such $R^*$ and $L^*$ are computed by the relation
\begin{equation*}
  L^*R=R^*L.
\end{equation*}
Each element above gives
\begin{gather}
  q_n-e_{n-1}^*=q_n^*-e_{n},\quad
  -q_ne_n^*=-q_{n+1}^*e_n\label{eq:drToda}
\end{gather}
for $n=0, 1, \dots, N-1$, where $e_{-1}=e_{-1}^*=e_{N-1}=e_{N-1}^*=0$.
Let us introduce
\begin{equation*}
  f_n\coloneq q_n+e_n,\quad n=0, 1, \dots, N-2,\quad
  f_{N-1}\coloneq q_{N-1}.
\end{equation*}
Then the relations \eqref{eq:drToda} yield
\begin{align*}
  f_n
  =q_n^*+e_{n-1}^*
  =q_n^*+e_{n-1}\frac{q_n^*}{q_{n-1}}
  =\frac{q_n^*}{q_{n-1}}(q_{n-1}+e_{n-1})
  =f_{n-1}\frac{q_n^*}{q_{n-1}}.
\end{align*}
Hence, we obtain
\begin{gather*}
  f_n=q_n+e_n,\quad n=0, 1, \dots, N-2,\quad f_{N-1}=q_{N-1},\\
  q_0^*=f_0,\quad q_n^*=q_{n-1}\frac{f_n}{f_{n-1}},\quad n=1, 2, \dots, N-1,\\
  e_n^*=e_n\frac{f_{n+1}}{f_n},\quad n=0, 1, \dots, N-2.
\end{gather*}

Let us iterate the procedure above.
Set $p^{(0)}_n(x)\coloneq p_n(x)$, $q_n^{(0)}\coloneq q_n$ and $e_n^{(0)}\coloneq e_n$.
Generate new monic polynomials $\{p_n^{(k)}(x)\}_{n=0}^N$, $k=1, 2, 3, \dots$, by
\begin{gather}
  p_0^{(k+1)}(x)\coloneq 1,\quad p^{(k+1)}_N(x)\coloneq p_N(x),\quad
  p_n^{(k+1)}(x)\coloneq p_n^{(k)}(x)-e_{n-1}^{(k)}p_{n-1}^{(k)}(x),\quad
  n=1, 2, \dots, N-1,\label{eq:GT-k-LBP}
\end{gather}
where $\{q_n^{(k)}\}_{n=0}^{N-1}$ and $\{e_n^{(k)}\}_{n=0}^{N-2}$ are computed by
\begin{gather}
  f_n^{(k)}=q_n^{(k)}+e_n^{(k)},\quad n=0, 1, \dots, N-2,\quad f_{N-1}^{(k)}=q_{N-1}^{(k)},\label{eq:drToda-d}\\
  q_0^{(k+1)}=f_0^{(k)},\quad q_n^{(k+1)}=q_{n-1}^{(k)}\frac{f_n^{(k)}}{f_{n-1}^{(k)}},\quad n=1, 2, \dots, N-1,\\
  e_n^{(k+1)}=e_n^{(k)}\frac{f_{n+1}^{(k)}}{f_n^{(k)}},\quad n=0, 1, \dots, N-2.\label{eq:drToda-e}
\end{gather}
The equations~\eqref{eq:drToda-d}--\eqref{eq:drToda-e} are called
the \emph{discrete relativistic Toda lattice}.
One can see that this system is a variation of the differential qd (dqd) algorithm~\cite{fernando1994asv,rutishauser1990lnm},
which is the subtraction-free version of the quotient-difference algorithm used to
compute the eigenvalues of a tridiagonal matrix.

Then $\{p_n^{(k)}(x)\}_{n=0}^N$ satisfy the relation
\begin{equation}\label{eq:CT-k-LBP}
  x p_n^{(k+1)}(x)=p_{n+1}^{(k)}(x)+q_n^{(k)}p_n^{(k)}(x),\quad
  n=0, 1, \dots, N-1,
\end{equation}
and the three-term recurrence relation
\begin{equation*}
  p^{(k)}_{n+1}(x)=(x-q_n^{(k)})p_n^{(k)}(x)-x e_{n-1}^{(k)} p_{n-1}^{(k)}(x),\quad
  n=0, 1, \dots, N-1.
\end{equation*}

\subsection{Isospectral transformation to a tridiagonal matrix}
Subtraction of \eqref{eq:GT-k-LBP} with $n\to n+1$ from \eqref{eq:CT-k-LBP} and $k\to k+n-1$ yield
\begin{equation}\label{eq:CT-k-OPS}
  x p_{n}^{(k+n)}(x)=p_{n+1}^{(k+n)}(x)+f_n^{(k+n-1)} p_n^{(k+n-1)}(x),\quad
  n=0, 1, \dots, N-1.
\end{equation}
Let us introduce variables
\begin{equation*}
  \hat q_n^{(k)}\coloneq f_n^{(k+n-1)},\quad
  \hat e_n^{(k)}\coloneq e_n^{(k+n)},
\end{equation*}
bidiagonal matrices
\begin{equation}\label{eq:bidiagonal-OPS}
  \hat R^{(k)}\coloneq
  \begin{pmatrix}
    \hat q_0^{(k)} & 1\\
    & \hat q_1^{(k)} & 1\\
    && \ddots & \ddots\\
    &&& \ddots & 1\\
    &&&& \hat q_{N-1}^{(k)}
  \end{pmatrix},\quad
  \hat L^{(k)}\coloneq
  \begin{pmatrix}
    1\\
    \hat e_0^{(k)} & 1\\
    & \hat e_1^{(k)} & \ddots\\
    && \ddots & \ddots \\
    &&& \hat e_{N-2}^{(k)} & 1
  \end{pmatrix},
\end{equation}
and vectors of polynomials
\begin{equation}\label{eq:vector-OPS}
  \hat{\bm p}^{(k)}(x)\coloneq
  \begin{pmatrix}
    \hat p_0^{(k)}(x)\\
    \hat p_1^{(k)}(x)\\
    \vdots\\
    \hat p_{N-1}^{(k)}(x)
  \end{pmatrix},
\end{equation}
where
\begin{equation*}
  \hat p_n^{(k)}(x)\coloneq p_n^{(k+n-1)}(x).
\end{equation*}
Then the relations~\eqref{eq:GT-k-LBP} and \eqref{eq:CT-k-OPS} yield
\begin{equation}\label{eq:GT-CT-OPS}
  \hat{\bm p}^{(k-1)}(x)=\hat L^{(k-1)}\hat{\bm p}^{(k)}(x),\quad
  x\hat{\bm p}^{(k+1)}(x)=\hat R^{(k)}\hat{\bm p}^{(k)}(x)+\bm p_N(x).
\end{equation}
Hence, $\hat{\bm p}^{(k)}(x)$ satisfies
\begin{equation}\label{eq:trr-OPS-vec}
  \hat T^{(k)}\hat{\bm p}^{(k)}(x)+\bm p_N(x)=x\hat{\bm p}^{(k)}(x),
\end{equation}
where $\hat T^{(k)}$ is a tridiagonal matrix
\begin{equation*}
  \hat T^{(k)}\coloneq \hat L^{(k)}\hat R^{(k)}=\hat R^{(k-1)}\hat L^{(k-1)}.
\end{equation*}
Therefore, Favard's theorem says that there exists a linear functional
$\mathcal L^{(k)}\colon \mathbb C[x]\to \mathbb C$ satisfying
the orthogonality condition
\begin{equation*}
  \mathcal L^{(k)}[\hat p_m^{(k)}(x)\hat p_n^{(k)}(x)]=\hat h_n^{(k)} \delta_{m, n},\quad
  \hat h_n^{(k)}\ne 0,\quad m, n=0, 1, \dots, N-1,
\end{equation*}
where $\delta_{m, n}$ is the Kronecker delta, and the terminating condition
\begin{equation*}
  \mathcal L^{(k)}[p_N(x)\pi(x)]=0 \quad \text{for all $\pi(x)\in\mathbb C[x]$}.
\end{equation*}
Note that the orthogonality relation is equivalent to
\begin{equation}\label{eq:orthogonality-OPS}
  \mathcal L^{(k)}[x^m \hat p_n^{(k)}(x)]=h_n^{(k)}\delta_{m, n},\quad
  h_n^{(k)}\ne 0,\quad
  n=0, 1, \dots, N-1,\quad m=0, 1, \dots, n.
\end{equation}
The linear functional $\mathcal L^{(k)}$ is concretely given by~\cite{akhiezer1965cmp}
\begin{equation}\label{eq:def-lf-tridiagonal-OPS}
  \mathcal L^{(k)}[\pi(x)]=\bm e_0^{\mathrm T}\pi(\hat T^{(k)})\bm e_0\quad
  \text{for all $\pi(x)\in\mathbb C[x]$},
\end{equation}
where
\begin{equation*}
  \bm e_0\coloneq
  \begin{pmatrix}
    1\\
    0\\
    \vdots\\
    0
  \end{pmatrix}\in \mathbb C^N.
\end{equation*}
Hence, if the eigenvalues $x_0, x_1, \dots, x_{N-1}$ of $\hat T^{(k)}$,
which are the zeros of $p_N(x)$,
are all simple, then there exist some constants $w_0^{(k)}, w_1^{(k)}, \dots, w_{N-1}^{(k)} \in \mathbb C$
such that
\begin{equation*}
  \mathcal L^{(k)}[\pi(x)]=\sum_{i=0}^{N-1} \pi(x_i)w_i^{(k)}\quad \text{for all $\pi(x)\in\mathbb C[x]$}.
\end{equation*}
This means that $\{\hat p_n^{(k)}(x)\}_{n=0}^N$ is
the \emph{monic discrete orthogonal polynomial sequence} with respect to $\mathcal L^{(k)}$.

In the theory of orthogonal polynomials, the relations~\eqref{eq:CT-k-OPS} and \eqref{eq:GT-k-LBP}
are called \emph{Christoffel transformation} and \emph{Geronimus transformation}, respectively~\cite{spiridonov1995ddt},
and the relation
\begin{equation*}
  \mathcal L^{(k+1)}[\pi(x)]=\mathcal L^{(k)}[x\pi(x)]\quad \text{for all $\pi(x)\in\mathbb C[x]$,}
\end{equation*}
is shown. Let us introduce the moment of the linear functional $\mathcal L^{(k)}$ as
\begin{equation}\label{eq:def-moment}
  \mu_{k+m}\coloneq \mathcal L^{(k)}[x^m]=\mathcal L^{(0)}[x^{k+m}].
\end{equation}
Then, from the orthogonality relation~\eqref{eq:orthogonality-OPS},
the determinant expression of the orthogonal polynomials $\{\hat p_n^{(k)}(x)\}_{n=0}^N$
is given by
\begin{equation*}
  \hat p_n^{(k)}(x)=\frac{1}{\tau^{(k)}_n}
  \begin{vmatrix}
    \mu_{k} & \mu_{k+1} & \dots & \mu_{k+n-1} & \mu_{k+n}\\
    \mu_{k+1} & \mu_{k+2} & \dots & \mu_{k+n} & \mu_{k+n+1}\\
    \vdots & \vdots & & \vdots & \vdots\\
    \mu_{k+n-1} & \mu_{k+n} & \dots & \mu_{k+2n-2} & \mu_{k+2n-1}\\
    1 & x & \dots & x^{n-1} & x^n
  \end{vmatrix}, \quad n=1, 2, \dots, N,
\end{equation*}
where $\tau_n^{(k)}$ is the Hankel determinant of the moments
\begin{equation*}
  \tau_0^{(k)}\coloneq 1,\quad
  \tau_n^{(k)}\coloneq |\mu_{k+i+j}|_{i, j=0}^{n-1},\quad
  n=1, 2, 3, \dots.
\end{equation*}

Next, we consider the determinant expression of $\{p_n^{(k)}(x)\}_{n=0}^N$:
\begin{equation*}
  p_n^{(k)}(x)=\hat p_n^{(k-n+1)}(x)=\frac{1}{\tau^{(k-n+1)}_n}
  \begin{vmatrix}
    \mu_{k-n+1} & \mu_{k-n+2} & \dots & \mu_{k} & \mu_{k+1}\\
    \mu_{k-n+2} & \mu_{k-n+3} & \dots & \mu_{k+1} & \mu_{k+2}\\
    \vdots & \vdots & & \vdots & \vdots\\
    \mu_{k} & \mu_{k+1} & \dots & \mu_{k+n-1} & \mu_{k+n}\\
    1 & x & \dots & x^{n-1} & x^n
  \end{vmatrix}, \quad n=1, 2, \dots, N.
\end{equation*}
Since we can extend the domain of the linear functional $\mathcal L^{(k)}$ for
$\mathbb C[x, x^{-1}]$ by \eqref{eq:def-lf-tridiagonal-OPS},
this determinant expression leads to the biorthogonality relation
\begin{equation*}
  \mathcal L^{(k)}[x^{-m} p_n^{(k)}(x)]=(-1)^n\frac{\tau_{n+1}^{(k-n)}}{\tau_n^{(k-n+1)}} \delta_{m, n},\quad
  n=0, 1, \dots, N-1,\quad m=0, 1, \dots, n.
\end{equation*}
The monic polynomials $\{p_n^{(k)}(x)\}_{n=0}^N$ are known as
the \emph{Laurent biorthogonal polynomials} with respect to $\mathcal L^{(k)}$.
Then, the relations \eqref{eq:CT-k-LBP}, \eqref{eq:GT-k-LBP}, \eqref{eq:CT-k-OPS} and
\begin{equation*}
  p_n^{(k)}(0)=(-1)^n \frac{\tau_n^{(k-n+2)}}{\tau_n^{(k-n+1)}}
\end{equation*}
give
\begin{gather*}
  q_n^{(k)}
  =-\frac{p_{n+1}^{(k)}(0)}{p_{n}^{(k)}(0)}
  =\frac{\tau_{n}^{(k-n+1)}\tau_{n+1}^{(k-n+1)}}{\tau_{n}^{(k-n+2)}\tau_{n+1}^{(k-n)}},\\
  e_n^{(k)}
  =-\frac{\mathcal L^{(k+1)}[x^{-n-1}p_{n+1}^{(k+1)}(x)]}{\mathcal L^{(k)}[x^{-n}p_{n}^{(k)}(x)]}
  =\frac{\tau_n^{(k-n+1)}\tau_{n+2}^{(k-n)}}{\tau_{n+1}^{(k-n)}\tau_{n+1}^{(k-n+1)}},\\
  f_n^{(k)}
  =-\frac{p_{n+1}^{(k+1)}(0)}{p_n^{(k)}(0)}
  =\frac{\tau_{n}^{(k-n+1)}\tau_{n+1}^{(k-n+2)}}{\tau_{n}^{(k-n+2)}\tau_{n+1}^{(k-n+1)}},
\end{gather*}
and
\begin{equation*}
  \hat q_n^{(k)}=f_n^{(k+n-1)}=\frac{\tau_{n}^{(k)}\tau_{n+1}^{(k+1)}}{\tau_{n}^{(k+1)}\tau_{n+1}^{(k)}},\quad
  \hat e_n^{(k)}=e_n^{(k+n)}=\frac{\tau_n^{(k+1)}\tau_{n+2}^{(k)}}{\tau_{n+1}^{(k)}\tau_{n+1}^{(k+1)}}.
\end{equation*}

Summarizing the above, we obtain Algorithm~\ref{alg:bi-bi-to-tri} of the isospectral transformation.

\begin{algorithm}[t]
  \caption{Isospectral transformation from bidiagonal--bidiagonal matrix pencil to tridiagonal matrix}
  \label{alg:bi-bi-to-tri}
\begin{algorithmic}
  \REQUIRE $\{q_n^{(0)}\}_{n=0}^{N-1}$ and $\{e_n^{(0)}\}_{n=0}^{N-2}$ (or bidiagonal matrices $R^{(0)}$ and $L^{(0)}$)
  \FOR{$k=0$ to $N-1$}
  \FOR{$n=0$ to $N-1$}
  \STATE\IfThenElse{$n<N-1$}
                   {$f_n^{(k)} \leftarrow q_n^{(k)}+e_n^{(k)}$}
                   {$f_n^{(k)} \leftarrow q_n^{(k)}$}
  \ENDFOR
  \FOR{$n=0$ to $N-1$}
  \STATE\IfThenElse{$n=0$}
                   {$q_n^{(k+1)} \leftarrow f_n^{(k)}$}
                   {$q_n^{(k+1)} \leftarrow q_{n-1}^{(k)}\frac{f_n^{(k)}}{f_{n-1}^{(k)}}$}
  \STATE\IfThen{$n<N-1$}
               {$e_n^{(k+1)} \leftarrow e_n^{(k)}\frac{f_{n+1}^{(k)}}{f_n^{(k)}}$}
  \ENDFOR
  \ENDFOR
  \FOR{$n=0$ to $N-1$}
  \STATE $\hat q_n^{(1)} \leftarrow f_n^{(n)}$
  \STATE\IfThen{$n<N-1$}
               {$\hat e_n^{(1)} \leftarrow e_n^{(n+1)}$}
  \ENDFOR
  \ENSURE $\{\hat q_n^{(1)}\}_{n=0}^{N-1}$ and $\{\hat e_n^{(1)}\}_{n=0}^{N-2}$ (or tridiagonal matrix $\hat T^{(1)}=\hat L^{(1)}\hat R^{(1)}$)
\end{algorithmic}
\end{algorithm}

\subsection{Numerical example}
Let us consider
\begin{equation*}
  R^{(0)}=
  \begin{pmatrix}
    1 & 1 & 0 & 0 & 0\\
    0 & 2 & 1 & 0 & 0\\
    0 & 0 & 3 & 1 & 0\\
    0 & 0 & 0 & 4 & 1\\
    0 & 0 & 0 & 0 & 5
  \end{pmatrix},\quad
  L^{(0)}=
  \begin{pmatrix}
    1 & 0 & 0 & 0 & 0\\
    -6 & 1 & 0 & 0 & 0\\
    0 & -7 & 1 & 0 & 0\\
    0 & 0 & -8 & 1 & 0\\
    0 & 0 & 0 & -9 & 1
  \end{pmatrix}.
\end{equation*}
Then, Algorithm~\ref{alg:bi-bi-to-tri} yields
\begin{gather*}
  \hat R^{(1)}=
  \begin{pmatrix}
    7& 1& 0& 0& 0\\
    0& \frac{620}{63}& 1& 0& 0\\
    0& 0& \frac{41949}{5890}& 1& 0\\
    0& 0& 0& \frac{5722439}{7639379}& 1\\
    0& 0& 0& 0& \frac{98340}{301181}
  \end{pmatrix},\quad
  \hat L^{(1)}=
  \begin{pmatrix}
    1& 0& 0& 0& 0\\
    \frac{54}{7}& 1& 0& 0& 0\\
    0& \frac{931}{90}& 1& 0& 0\\
    0& 0& \frac{4720320}{2745329}& 1& 0\\
    0& 0& 0& \frac{90306875}{493635659}& 1
  \end{pmatrix},
\end{gather*}
and
\begin{gather*}
  \hat T^{(1)}
  =\hat L^{(1)}\hat R^{(1)}
  =
  \begin{pmatrix}
    7& 1& 0& 0& 0\\[.5em]
    54& \frac{158}{9}& 1& 0& 0\\[.5em]
    0& \frac{8246}{81}& \frac{92590}{5301}& 1& 0\\[.5em]
    0& 0& \frac{4248288}{346921}& \frac{2382991}{965371} & 1\\[.5em]
    0& 0& 0& \frac{368125}{2686321} & \frac{835}{1639}
  \end{pmatrix}
\end{gather*}
The computation above was executed by SymPy (a symbolic computation library for Python).
The eigenvalues of both $(R^{(0)}, L^{(0)})$ and $\hat T^{(1)}$ are
29.10515103, 12.22484344, 2.82190399, 0.17848385 and 0.66961769,
which are verified by \texttt{numpy.linalg.eig()}.

\section{Isospectral transformation between tridiagonal matrix and tridiagonal--bidiagonal matrix pencil}
\label{sec:transf-betw-trid-1}

\subsection{Sequence of tridiagonal--bidiagonal GEVPs and eigenvectors}
We first fix a vector
$\bm{\varepsilon} = (\varepsilon_0,\varepsilon_1,\ldots,\varepsilon_{N-2}) \in \{0,1\}^{N-1}$.
Let $\{p_n^{(k)}(x)\}_{n=0}^N$ be monic polynomials of degree $n$ defined by
the three-term recurrence relation
\begin{gather}
  p_{-1}^{(k)}(x)\coloneq 0,\quad p_0^{(k)}(x)\coloneq 1,\nonumber\\
  \begin{multlined}[b]
    p_{n+1}^{(k)}(x)\coloneq (x-q_n^{(k)}-(1-\epsilon_{n-1})e_{n-1}^{(k)})p_n^{(k)}(x)-(\epsilon_{n-1}x+(1-\epsilon_{n-1})q_{n-1}^{(k)})e_{n-1}^{(k)}p_{n-1}^{(k)}(x),\\
    n=0, 1, \dots, N-1,
  \end{multlined}\label{eq:trr-e-BLP}
\end{gather}
where $e_{-1}\coloneq 0$.
Let us consider lower bidiagonal matrices
\begin{equation}\label{eq:def-L-epsilon}
  L_{(\xi_0,\ldots,\xi_{N-2})}^{(k)}\coloneq
  \begin{pmatrix}
    1\\
    -\xi_0 e_0^{(k)} & 1\\
    & -\xi_1 e_1^{(k)} & \ddots\\
    && \ddots & \ddots\\
    &&& -\xi_{N-2} e_{N-2}^{(k)} & 1
  \end{pmatrix},
\end{equation}
upper bidiagonal matrices
\begin{equation*}
  R^{(k)}\coloneq
  \begin{pmatrix}
    q_0^{(k)} & 1\\
    & q_1^{(k)} & 1\\
    && \ddots & \ddots\\
    &&& \ddots & 1\\
    &&&& q_{N-1}^{(k)}
  \end{pmatrix},
\end{equation*}
and vectors of polynomials which are the same as \eqref{eq:def-vec-p}.
Then we can rewrite the three-term recurrence relation~\eqref{eq:trr-e-BLP} as
\begin{equation}\label{eq:trr-e-BLP-vec}
  L_{\bm{\varepsilon}^*}^{(k)}R^{(k)}\bm p^{(k)}(x)+\bm p_N(x)=x L_{\bm\varepsilon}^{(k)}\bm p^{(k)}(x),
\end{equation}
where $\bm{\varepsilon}^*\coloneq\bm{\varepsilon}-(1,1,\ldots,1)$.

\begin{example} When $N=5, \bm\varepsilon=(1,1,0,1)$,
  we obtain
  \begin{equation*}
    L_{\bm{\varepsilon}^*}^{(k)} =
    \begin{pmatrix}
      1 \\
      0 & 1\\
      & 0 & 1\\\
      && e_2^{(k)} & 1\\
      &&& 0 & 1
    \end{pmatrix}, \quad
    L_{\bm{\varepsilon}}^{(k)} =
    \begin{pmatrix}
      1 \\
      -e_0^{(k)} & 1\\
      & -e_1^{(k)} & 1\\\
      && 0 & 1\\
      &&& -e_3^{(k)} & 1
    \end{pmatrix}.
  \end{equation*}
\end{example}

Note that
\begin{itemize}
\item if $\bm\varepsilon=(0, 0, \dots, 0)$, then \eqref{eq:trr-e-BLP-vec} becomes
  the three-term recurrence relation of monic orthogonal polynomials (see \eqref{eq:trr-OPS-vec});
\item if $\bm\varepsilon=(1, 1, \dots, 1)$, then \eqref{eq:trr-e-BLP-vec} becomes
  the three-term recurrence relation of monic Laurent biorthogonal polynomials (see \eqref{eq:trr-LBP-vec}).
\end{itemize}
Hence the monic polynomial sequence $\{p_n^{(k)}(x)\}_{n=0}^N$ is a generalization
of both the orthogonal polynomials and the Laurent biorthogonal polynomials.
We therefore use the notation that is almost the same as that of the previous section.

Let us introduce new monic polynomials $\{p_n^{(k+1)}(x)\}_{n=0}^N$ of degree $n$
generated from the monic polynomials $\{p_n^{(k)}(x)\}_{n=0}^N$ by
\begin{equation*}
  \bm p^{(k+1)}(x)\coloneq \left(L_{\bm{\varepsilon}^*}^{(k)}\right)^{-1}L_{\bm{\varepsilon}}^{(k)}\bm p^{(k)}(x).
\end{equation*}
Substituting $\bm p^{(k)}(x)=\left(L_{\bm{\varepsilon}}^{(k)}\right)^{-1}L_{\bm{\varepsilon}^*}^{(k)}\bm p^{(k+1)}(x)$ into
\eqref{eq:trr-e-BLP-vec}, we obtain
\begin{equation}\label{eq:GT-CT-BLP}
  L_{\bm{\varepsilon}^*}^{(k)}\bm p^{(k+1)}(x)=L_{\bm\varepsilon}^{(k)}\bm p^{(k)}(x),\quad
  R^{(k)}\bm p^{(k)}(x)+\bm p_N(x)=x\bm p^{(k+1)}(x).
\end{equation}
This is a generalization of \eqref{eq:GT-CT-LBP} and \eqref{eq:GT-CT-OPS}.
Elimination of $\bm p(x)$ yields
\begin{equation*}
  R^{(k)}\left(L_{\bm{\varepsilon}}^{(k)}\right)^{-1}L_{\bm{\varepsilon}^*}^{(k)}\bm p^{(k+1)}(x)+\bm p_N(x)=x\bm p^{(k+1)}(x).
\end{equation*}
The LU factorization generates bidiagonal matrices
$L_{\bm{\varepsilon}}^{(k+1)}$, $L_{\bm{\varepsilon}^*}^{(k+1)}$, $\tilde R^{(k)}$ and $R^{(k+1)}$ satisfying
\begin{equation}\label{eq:de-Toda}
  L_{\bm{\varepsilon}}^{(k+1)}R^{(k)}=\tilde R^{(k)} L_{\bm{\varepsilon}}^{(k)},\quad
  \tilde R^{(k)} L_{\bm{\varepsilon}^*}^{(k)}=L_{\bm{\varepsilon}^*}^{(k+1)}R^{(k+1)}.
\end{equation}
Then, we obtain
\begin{equation*}
  L_{\bm{\varepsilon}^*}^{(k+1)}R^{(k+1)}\bm p^{(k+1)}(x)+\bm p_N(x)=x L_{\bm{\varepsilon}}^{(k+1)}\bm p^{(k+1)}(x).
\end{equation*}
Each element of \eqref{eq:de-Toda} give the relations
\begin{gather*}
  q_n^{(k)}-\epsilon_{n-1}e_{n-1}^{(k+1)}=\tilde q_n^{(k)}-\epsilon_n e_n^{(k)},\quad
  \tilde q_n^{(k)}+(1-\epsilon_n) e_n^{(k)}=q_n^{(k+1)}+(1-\epsilon_{n-1})e_{n-1}^{(k+1)}
\end{gather*}
for $n=0, 1, \dots, N-1$, and
\begin{gather*}
  q_n^{(k)} e_n^{(k+1)}=\tilde q_{n+1}^{(k)} e_n^{(k)} \quad \text{if $\epsilon_n=1$},\\
  \tilde q_{n+1}^{(k)} e_n^{(k)}=q_n^{(k+1)} e_n^{(k+1)} \quad \text{if $\epsilon_n=0$}
\end{gather*}
for $n=0, 1, \dots, N-2$.
Elimination of $\tilde q_n^{(k)}$ yields
\begin{gather}
  q_n^{(k)}+e_n^{(k)}=q_n^{(k+1)}+e_{n-1}^{(k+1)},\label{eq:de-Toda-p}\\
  q_n^{(k)}e_n^{(k+1)}=(q_{n+1}^{(k+1)}-(1-\epsilon_{n+1})e_{n+1}^{(k)})e_n^{(k)}\quad \text{if $\epsilon_n=1$}\label{eq:de-Toda-m1},\\
  (q_{n+1}^{(k)}+\epsilon_{n+1}e_{n+1}^{(k)})e_n^{(k)}=q_n^{(k+1)}e_n^{(k+1)}\quad \text{if $\epsilon_n=0$}.\label{eq:de-Toda-m2}
\end{gather}
Further, substituting \eqref{eq:de-Toda-p} to \eqref{eq:de-Toda-m1} and \eqref{eq:de-Toda-m2},
we can see that
\begin{gather}
  (q_n^{(k)}-(1-\epsilon_n)e_{n-1}^{(k+1)})e_n^{(k+1)}=(q_{n+1}^{(k+1)}-(1-\epsilon_{n+1})e_{n+1}^{(k)})e_n^{(k)},\label{eq:de-Toda-m3}\\
  (q_n^{(k+1)}+\epsilon_{n}e_{n-1}^{(k+1)})e_n^{(k+1)}=(q_{n+1}^{(k)}+\epsilon_{n+1}e_{n+1}^{(k)})e_n^{(k)}\label{eq:de-Toda-m4}
\end{gather}
hold in any case.
Let us introduce
\begin{gather*}
  d_0^{(k+1)}\coloneq q_0^{(k)}+\epsilon_0 e_0^{(k)},\quad
  d_n^{(k+1)}\coloneq q_n^{(k+1)}-(1-\epsilon_n) e_{n}^{(k)},\quad n=1, \dots, N-1,\\
  f_n^{(k)}\coloneq q_n^{(k)}+\epsilon_n e_n^{(k)},\quad n=0, 1, \dots, N-2,\quad
  f_{N-1}^{(k)}\coloneq q_{N-1}^{(k)}.
\end{gather*}
Then, from \eqref{eq:de-Toda-p}, \eqref{eq:de-Toda-m3} and \eqref{eq:de-Toda-m4}, we obtain
\begin{align*}
  d_n^{(k+1)}
  &=q_n^{(k)}+\epsilon_n e_n^{(k)}-e_{n-1}^{(k+1)}\\
  &=\frac{f_n^{(k)}}{q_{n-1}^{(k+1)}+\epsilon_{n-1}e_{n-2}^{(k+1)}}(q_{n-1}^{(k+1)}+\epsilon_{n-1}e_{n-2}^{(k+1)}-e_{n-1}^{(k)})\\
  &=\frac{f_n^{(k)}}{q_{n-1}^{(k+1)}+\epsilon_{n-1}e_{n-2}^{(k+1)}}(d_{n-1}^{(k+1)}+\epsilon_{n-1}(e_{n-2}^{(k+1)}-e_{n-1}^{(k)})),\\
  f_n^{(k)}
  &=q_n^{(k+1)}-(1-\epsilon_n)e_n^{(k)}+e_{n-1}^{(k+1)}\\
  &=\frac{d_n^{(k+1)}}{q_{n-1}^{(k)}-(1-\epsilon_{n-1})e_{n-2}^{(k+1)}}(q_{n-1}^{(k)}-(1-\epsilon_{n-1})e_{n-2}^{(k+1)}+e_{n-1}^{(k)})\\
  &=\frac{d_n^{(k+1)}}{q_{n-1}^{(k)}-(1-\epsilon_{n-1})e_{n-2}^{(k+1)}}(f_{n-1}^{(k)}-(1-\epsilon_{n-1})(e_{n-2}^{(k+1)}-e_{n-1}^{(k)})).
\end{align*}
Hence, we can compute $\{q_n^{(k+1)}\}_{n=0}^{N-1}$ and $\{e_n^{(k+1)}\}_{n=0}^{N-2}$
from $\{q_n^{(k)}\}_{n=0}^{N-1}$ and $\{e_n^{(k)}\}_{n=0}^{N-2}$ by
\begin{gather*}
  f_n^{(k)}=q_n^{(k)}+\epsilon_n e_n^{(k)},\quad n=0, 1, \dots, N-2,\quad
  f_{N-1}^{(k)}=q_{N-1}^{(k)},\\
  d_0^{(k+1)}=f_0^{(k)},\quad
  d_n^{(k+1)}=
  \begin{dcases*}
    d_{n-1}^{(k+1)}\frac{f_n^{(k)}}{q_{n-1}^{(k+1)}} & if $\epsilon_{n-1}=0$,\\
    q_{n-1}^{(k)}\frac{f_n^{(k)}}{f_{n-1}^{(k)}} & if $\epsilon_{n-1}=1$,
  \end{dcases*}\quad
  n=1, 2, \dots, N-1,\\
  q_n^{(k+1)}=d_n^{(k+1)}+(1-\epsilon_n)e_n^{(k)},\quad n=0, 1, \dots, N-1,\\
  e_n^{(k+1)}=e_n^{(k)}\frac{f_{n+1}^{(k)}}{q_n^{(k+1)}+\epsilon_n e_{n-1}^{(k+1)}},\quad
  n=0, 1, \dots, N-2.
\end{gather*}
We call this system the \emph{discrete elementary Toda orbits}~\cite{kobayashi2021nde}.

\subsection{Isospectral transformation to a tridiagonal matrix}
Each element of \eqref{eq:GT-CT-BLP} gives
\begin{alignat}{2}
  xp_{n}^{(k+1)}(x)&=p_{n+1}^{(k)}(x)+q_n^{(k)}p_n^{(k)}(x),\label{eq:e-BLP-st-q}\\
  p_{n+1}^{(k-1)}(x)&=p_{n+1}^{(k)}(x)+e_n^{(k-1)}p_n^{(k)}(x) &\quad& \text{if $\epsilon_n=0$},\label{eq:e-BLP-st-e0}\\
  p_{n+1}^{(k+1)}(x)&=p_{n+1}^{(k)}(x)-e_n^{(k)}p_n^{(k)}(x) &\quad& \text{if $\epsilon_n=1$}.\label{eq:e-BLP-st-e1}\\
\intertext{Therefore, we obtain}
  xp_{n}^{(k+1)}(x)&=p_{n+1}^{(k-1)}(x)+d_n^{(k)}p_n^{(k)}(x) &\quad& \text{if $\epsilon_n=0$},\label{eq:e-BLP-st-d0}\\
  xp_{n}^{(k+1)}(x)&=p_{n+1}^{(k+1)}(x)+f_n^{(k)}p_n^{(k)}(x) && \text{if $\epsilon_n=1$}.\label{eq:e-BLP-st-f1}
\end{alignat}
Since $d_n^{(k)}=q_n^{(k)}$ if $\epsilon_n=1$ and $f_n^{(k)}=q_n^{(k)}$ if $\epsilon_n=0$,
we can see that
\begin{gather*}
  xp_n^{(k+1)}(x)=p_{n+1}^{(k+\epsilon_n-1)}(x)+d_n^{(k)}p_n^{(k)}(x),\\
  xp_n^{(k+1)}(x)=p_{n+1}^{(k+\epsilon_n)}(x)+f_n^{(k)}p_n^{(k)}(x)
\end{gather*}
hold in any case.
It also holds in any case that
\begin{equation*}
  p_n^{(k)}(x)=p_n^{(k+1)}(x)+e_{n-1}^{(k)} p_{n-1}^{(k+1-\epsilon_{n-1})}(x).
\end{equation*}
These relations lead to the three-term recurrence relation
\begin{equation}\label{eq:trr-e-BLP-to-OPS-pre}
  xp_n^{(k)}(x)=p_{n+1}^{(k+\epsilon_n)}(x)+(f_n^{(k)}+e_{n-1}^{(k)})p_n^{(k)}(x)+f_{n-1}^{(k-\epsilon_{n-1})}e_{n-1}^{(k)}p_{n-1}^{(k-\epsilon_{n-1})}(x).
\end{equation}

Let us introduce $\bm{\eta}=(\eta_0,\eta_1,\ldots,\eta_{N-1})$ by
\begin{equation*}
  \eta_0\coloneq 0,\quad
  \eta_n\coloneq \sum_{j=0}^{n-1}\epsilon_j,\quad
  n=1, 2, \dots, N-1,
\end{equation*}
variables
\begin{equation*}
  \hat q_n^{(k)}\coloneq f_n^{(k+\eta_n)},\quad
  \hat e_n^{(k)}\coloneq e_n^{(k+\eta_{n+1})},
\end{equation*}
and polynomials
\begin{equation*}
  \hat p_n^{(k)}(x)\coloneq p_n^{(k+\eta_n)}(x).
\end{equation*}
Then the relation \eqref{eq:trr-e-BLP-to-OPS-pre} yields
\begin{equation*}
  \hat T^{(k)}\hat{\bm p}^{(k)}(x)+\bm p_N(x)=x\hat{\bm p}^{(k)}(x),
\end{equation*}
that is the same as \eqref{eq:trr-OPS-vec};
the monic polynomial sequence $\{\hat p_n^{(k)}(x)\}_{n=0}^N$ is
the orthogonal polynomial sequence with respect to the linear functional
$\mathcal L^{(k)}$ defined by \eqref{eq:def-lf-tridiagonal-OPS}.
This fact gives the determinant expression of $\{p_n^{(k)}(x)\}_{n=0}^N$:
\begin{equation*}
  p_n^{(k)}(x)=\frac{1}{\tau_n^{(k-\eta_n)}}
  \begin{vmatrix}
    \mu_{k-\eta_n} & \mu_{k-\eta_n+1} & \dots & \mu_{k-\eta_n+n-1} & \mu_{k-\eta_n+n}\\
    \mu_{k-\eta_n+1} & \mu_{k-\eta_n+2} & \dots & \mu_{k-\eta_n+n} & \mu_{k-\eta_n+n+1}\\
    \vdots & \vdots & & \vdots & \vdots\\
    \mu_{k-\eta_n+n-1} & \mu_{k-\eta_n+n} & \dots & \mu_{k-\eta_n+2n-2} & \mu_{k-\eta_n+2n-1}\\
    1 & x & \dots & x^{n-1} & x^n
  \end{vmatrix},\quad
  n=1, 2, \dots, N.
\end{equation*}
Hence, $\{p_n^{(k)}(x)\}_{n=0}^N$ satisfies the $\epsilon$-biorthogonality relation
for $n=0, 1, \dots, N-1$:
\begin{equation*}
  \mathcal L^{(k)}[x^m p_n^{(k)}(x)]=0,\quad  m=-\eta_n, -\eta_n+1, \dots, -\eta_n+n-1,
\end{equation*}
and
\begin{equation*}
  \mathcal L^{(k)}[x^{-\eta_n+n-\epsilon_n(n+1)} p_n^{(k)}(x)]=(1-2\epsilon_n)^n\frac{\tau_{n+1}^{(k-\eta_{n+1})}}{\tau_n^{(k-\eta_n)}},
\end{equation*}
or
\begin{gather*}
  \mathcal L^{(k)}[x^{-\eta_n+n} p_n^{(k)}(x)]=\frac{\tau_{n+1}^{(k-\eta_{n})}}{\tau_n^{(k-\eta_n)}},\quad
  \mathcal L^{(k)}[x^{-\eta_n-1} p_n^{(k)}(x)]=(-1)^n\frac{\tau_{n+1}^{(k-\eta_{n}-1)}}{\tau_n^{(k-\eta_n)}}.
\end{gather*}
The monic polynomial sequence $\{p_n^{(k)}(x)\}_{n=0}^N$ is called
the \emph{$\epsilon$-biorthogonal Laurent polynomials}~\cite{faybusovich2001imp}.
Then, the relations \eqref{eq:e-BLP-st-q}--\eqref{eq:e-BLP-st-f1} and
\begin{equation*}
  p_n^{(k)}(0)=(-1)^n\frac{\tau_n^{(k-\eta_n+1)}}{\tau_n^{(k-\eta_n)}},\quad
  p_{n+1}^{(k+\epsilon_{n})}(0)=(-1)^{n+1}\frac{\tau_{n+1}^{(k-\eta_{n+1}+\epsilon_n+1)}}{\tau_{n+1}^{(k-\eta_{n+1}+\epsilon_n)}}=(-1)^{n+1}\frac{\tau_{n+1}^{(k-\eta_{n}+1)}}{\tau_{n+1}^{(k-\eta_{n})}}
\end{equation*}
give
\begin{gather*}
  q_n^{(k)}
  =-\frac{p_{n+1}^{(k)}(0)}{p_{n}^{(k)}(0)}
  =\frac{\tau_n^{(k-\eta_n)}\tau_{n+1}^{(k-\eta_{n+1}+1)}}{\tau_{n}^{(k-\eta_n+1)}\tau_{n+1}^{(k-\eta_{n+1})}},\\
  d_n^{(k)}
  =-\frac{p_{n+1}^{(k+\epsilon_n-1)}(0)}{p_{n}^{(k)}(0)}
  =\frac{\tau_n^{(k-\eta_n)}\tau_{n+1}^{(k-\eta_{n})}}{\tau_{n}^{(k-\eta_n+1)}\tau_{n+1}^{(k-\eta_{n}-1)}},\\
  f_n^{(k)}
  =-\frac{p_{n+1}^{(k+\epsilon_n)}(0)}{p_{n}^{(k)}(0)}
  =\frac{\tau_n^{(k-\eta_n)}\tau_{n+1}^{(k-\eta_{n}+1)}}{\tau_{n}^{(k-\eta_n+1)}\tau_{n+1}^{(k-\eta_{n})}}
\end{gather*}
and, since $\mathcal L^{(k+1)}[x^{-\eta_{n+1}+n-\epsilon_n(n+1)}\pi(x)]=\mathcal L^{(k+1-\epsilon_n)}[x^{-\eta_{n}+n-\epsilon_n(n+1)}\pi(x)]=\mathcal L^{(k)}[x^{-\eta_{n+1}+n+1-\epsilon_n(n+1)}\pi(x)]$
for all $\pi(x) \in \mathbb C[x]$ and
\begin{gather*}
  \mathcal L^{(k)}[x^{-\eta_{n+1}+n+1-\epsilon_n(n+1)}p_{n+1}^{(k)}(x)]=
  \begin{dcases*}
    \frac{\tau_{n+2}^{(k-\eta_{n+1})}}{\tau_{n+1}^{(k-\eta_{n+1})}} & if $\epsilon_n=0$,\\
    0 & if $\epsilon_n=1$,
  \end{dcases*}\\
  \mathcal L^{(k+1)}[x^{-\eta_{n+1}+n-\epsilon_n(n+1)}p_{n+1}^{(k+1)}(x)]=
  \begin{dcases*}
    0 & if $\epsilon_n=0$,\\
    (-1)^{n+1}\frac{\tau_{n+2}^{(k-\eta_{n+1})}}{\tau_{n+1}^{(k-\eta_{n+1}+1)}} & if $\epsilon_n=1$,
  \end{dcases*}\\
  \mathcal L^{(k+1-\epsilon_n)}[x^{-\eta_n+n-\epsilon_n(n+1)}p_n^{(k+1-\epsilon_n)}(x)]=
  \begin{dcases*}
    \frac{\tau_{n+1}^{(k-\eta_n+1)}}{\tau_n^{(k-\eta_n+1)}}=\frac{\tau_{n+1}^{(k-\eta_{n+1}+1)}}{\tau_n^{(k-\eta_{n+1}+1)}} & if $\epsilon_n=0$,\\
    (-1)^n\frac{\tau_{n+1}^{(k-\eta_n-1)}}{\tau_n^{(k-\eta_n)}}=(-1)^n\frac{\tau_{n+1}^{(k-\eta_{n+1})}}{\tau_n^{(k-\eta_{n+1}+1)}} & if $\epsilon_n=1$
  \end{dcases*}
\end{gather*}
hold,
\begin{equation*}
  e_n^{(k)}
  =\frac{\mathcal L^{(k)}[x^{-\eta_{n+1}+n+1-\epsilon_n(n+1)}p_{n+1}^{(k)}(x)]-\mathcal L^{(k+1)}[x^{-\eta_{n+1}+n-\epsilon_n(n+1)}p_{n+1}^{(k+1)}(x)]}{\mathcal L^{(k+1-\epsilon_n)}[x^{-\eta_n+n-\epsilon_n(n+1)}p_n^{(k+1-\epsilon_n)}(x)]}
  =\frac{\tau_n^{(k-\eta_{n+1}+1)}\tau_{n+2}^{(k-\eta_{n+1})}}{\tau_{n+1}^{(k-\eta_{n+1}+1)}\tau_{n+1}^{(k-\eta_{n+1})}}.
\end{equation*}

Summarizing the above, we obtain Algorithm~\ref{alg:tri-bi-to-tri} of
the isospectral transformation.

\begin{algorithm}[t]
  \caption{Isospectral transformation from tridiagonal--bidiagonal matrix pencil to tridiagonal matrix}
  \label{alg:tri-bi-to-tri}
\begin{algorithmic}
  \REQUIRE $\{q_n^{(0)}\}_{n=0}^{N-1}$, $\{e_n^{(0)}\}_{n=0}^{N-2}$ and $\bm\epsilon\in\{0, 1\}^{N-1}$ (or bidiagonal matrices $R^{(0)}$, $L_{\bm\epsilon^*}^{(0)}$ and $L_{\bm\epsilon}^{(0)}$)
  \STATE $\eta_0 \leftarrow 0$
  \FOR{$n=1$ to $N-1$}
  \STATE $\eta_n \leftarrow \eta_{n-1}+\epsilon_{n-1}$
  \ENDFOR
  \FOR{$k=0$ to $\eta_{N-1}$}
  \FOR{$n=0$ to $N-1$}
  \STATE\IfThenElse{$n<N-1$}
                   {$f_n^{(k)} \leftarrow q_n^{(k)}+\epsilon_n e_n^{(k)}$}
                   {$f_n^{(k)} \leftarrow q_n^{(k)}$}
  \ENDFOR
  \FOR{$n=0$ to $N-1$}
  \STATE\IfThenElse{$n=0$}
                   {$d_n^{(k+1)}\leftarrow f_n^{(k)}$}
                   {\IfThenElse{$\epsilon_{n-1}=0$}
                               {$d_n^{(k+1)}\leftarrow d_{n-1}^{(k+1)}\frac{f_n^{(k)}}{q_{n-1}^{(k+1)}}$}
                               {$d_n^{(k+1)}\leftarrow q_{n-1}^{(k)}\frac{f_n^{(k)}}{f_{n-1}^{(k)}}$}
                   }
  \STATE $q^{(k+1)}_n \leftarrow d_n^{(k+1)}+(1-\epsilon_n)e_n^{(k)}$
  \STATE\IfThen{$n<N-1$}
               {$e_n^{(k+1)} \leftarrow e_n^{(k)}\frac{f_{n+1}^{(k)}}{q_n^{(k+1)}+\epsilon_n e_{n-1}^{(k+1)}}$}
  \ENDFOR
  \ENDFOR
  \FOR{$n=0$ to $N-1$}
  \STATE $\hat q_n^{(0)} \leftarrow f_n^{(\eta_n)}$
  \STATE\IfThen{$n<N-1$}
               {$\hat e_n^{(0)} \leftarrow e_n^{(\eta_{n+1})}$}
  \ENDFOR
  \ENSURE $\{\hat q_n^{(0)}\}_{n=0}^{N-1}$ and $\{\hat e_n^{(0)}\}_{n=0}^{N-2}$ (or tridiagonal matrix $\hat T^{(0)}=\hat L^{(0)}\hat R^{(0)}$)
\end{algorithmic}
\end{algorithm}

\subsection{Numerical example}
Let us consider
\begin{gather*}
  L_{\bm\epsilon^*}^{(0)}=
  \begin{pmatrix}
    1 & 0 & 0 & 0 & 0 & 0\\
    0 & 1 & 0 & 0 & 0 & 0\\
    0 & 0 & 1 & 0 & 0 & 0\\
    0 & 0 & 0 & 1 & 0 & 0\\
    0 & 0 & 0 & 10 & 1 & 0\\
    0 & 0 & 0 & 0 & 11 & 1
  \end{pmatrix},\quad
  R^{(0)}=
  \begin{pmatrix}
    1 & 1 & 0 & 0 & 0 & 0\\
    0 & 2 & 1 & 0 & 0 & 0\\
    0 & 0 & 3 & 1 & 0 & 0\\
    0 & 0 & 0 & 4 & 1 & 0\\
    0 & 0 & 0 & 0 & 5 & 1\\
    0 & 0 & 0 & 0 & 0 & 6
  \end{pmatrix},
\end{gather*}
i.e.,
\begin{gather*}
  L_{\bm\epsilon^*}^{(0)}R^{(0)}
  =
  \begin{pmatrix}
    1 & 1 & 0 & 0 & 0 & 0\\
    0 & 2 & 1 & 0 & 0 & 0\\
    0 & 0 & 3 & 1 & 0 & 0\\
    0 & 0 & 0 & 4 & 1 & 0\\
    0 & 0 & 0 & 40 & 15 & 1\\
    0 & 0 & 0 & 0 & 55 & 17
  \end{pmatrix},
\end{gather*}
and
\begin{gather*}
  L_{\bm\epsilon}^{(0)}=
  \begin{pmatrix}
    1 & 0 & 0 & 0 & 0 & 0\\
    -7 & 1 & 0 & 0 & 0 & 0\\
    0 & -8 & 1 & 0 & 0 & 0\\
    0 & 0 & -9 & 1 & 0 & 0\\
    0 & 0 & 0 & 0 & 1 & 0\\
    0 & 0 & 0 & 0 & 0 & 1
  \end{pmatrix}.
\end{gather*}
Then, Algorithm~\ref{alg:tri-bi-to-tri} yields
\begin{gather*}
  \hat R^{(0)}=
  \begin{pmatrix}
    8& 1& 0& 0& 0& 0\\
    0& \frac{217}{20}& 1& 0& 0& 0\\
    0& 0& \frac{13150}{1953}& 1& 0& 0\\
    0& 0& 0& \frac{3924423}{614105}& 1& 0\\
    0& 0& 0& 0& \frac{2596480772}{1515844721}& 1\\
    0& 0& 0& 0& 0 & \frac{156435}{1389979}
  \end{pmatrix},\quad
  \hat L^{(0)}=
  \begin{pmatrix}
    1& 0& 0& 0& 0 &0\\
    \frac{35}{4} & 1& 0& 0& 0 &0\\
    0& \frac{5184}{1085}& 1& 0& 0 &0\\
    0& 0& \frac{101339}{11835}& 1& 0 &0\\
    0& 0& 0& \frac{685706750}{67877983}& 1& 0\\
    0& 0& 0& 0& \frac{119912925}{14496090991} & 1
  \end{pmatrix},
\end{gather*}
and
\begin{gather*}
  \hat T^{(0)}=\hat L^{(0)}\hat R^{(0)}=
  \begin{pmatrix}
    8& 1& 0& 0& 0 &0\\[0.5em]
    70& \frac{98}{5}& 1& 0& 0&0\\[0.5em]
    0& \frac{1296}{25}& \frac{518}{45}& 1& 0&0\\[0.5em]
    0& 0& \frac{4670}{81}& \frac{62848}{4203}& 1&0\\[0.5em]
    0& 0& 0& \frac{14079150}{218089}& \frac{57542826}{4870343}&1\\[0.5em]
    0 & 0& 0& 0& \frac{1541100}{108764041}& \frac{1260}{10429}\\[0.5em]
  \end{pmatrix}.
\end{gather*}
The eigenvalues of both $(L_{\bm\epsilon^*}^{(0)}R^{(0)}, L_{\bm\epsilon}^{(0)})$ and $\hat T^{(0)}$ are 28.1051142, 22.50730913, 10.85981143,
4.18583949, 0.23568694 and 0.10623881,
which are verified by \texttt{numpy.linalg.eig()}.

\section{Isospectral transformation between Hessenberg matrix and Hessenberg--bidiagonal matrix pencil}
\label{sec:transf-betw-hess}

\subsection{Sequence of Hessenberg--bidiagonal GEVPs and eigenvectors}
As a generalization, let us consider a Hessenberg--bidiagonal matrix pencil,
i.e., the Hessenberg matrix of the form
\begin{equation*}
  H^{(k)}\coloneq L_{\bm{\varepsilon}^*}^{(k)}R^{(k+M-1)}R^{(k+M-2)}\dots R^{(k)},
\end{equation*}
where $M$ is a positive integer,
and the bidiagonal matrix $L_{\bm\varepsilon}^{(k)}$ defined by \eqref{eq:def-L-epsilon}.
To obtain the polynomial sequence related to the matrix pencil,
let us consider polynomial sequences $\{p_n^{(k)}(x)\}_{n=0}^N$ satisfying the relations
\begin{equation*}
  L_{\bm{\varepsilon}^*}^{(k)}\bm p^{(k+M)}(x)=L_{\bm\varepsilon}^{(k)}\bm p^{(k)}(x),\quad
  R^{(k)}\bm p^{(k)}(x)+\bm p_N^{(k)}(x)=x\bm p^{(k+1)}(x),
\end{equation*}
and
\begin{equation*}
  p_N^{(k+M)}(x)=p_N^{(k)}(x)
\end{equation*}
for all integer $k$.
Note that $p_N^{(k+1)}(x)=p_N^{(k)}(x)$ does not hold in general.
Then the matrix equation~\eqref{eq:de-Toda} is generalized as
\begin{equation*}
  L_{\bm{\varepsilon}}^{(k+1)}R^{(k)}=\tilde R^{(k)} L_{\bm{\varepsilon}}^{(k)},\quad
  \tilde R^{(k)} L_{\bm{\varepsilon}^*}^{(k)}=L_{\bm{\varepsilon}^*}^{(k+1)}R^{(k+M)},
\end{equation*}
and it is shown that $\{p_n^{(k)}(x)\}_{n=0}^N$ satisfies
the $(M+2)$-terms recurrence relation
\begin{equation*}
  H^{(k)}\bm p^{(k)}(x)+\sum_{j=0}^{M-1}L_{\bm{\varepsilon}^*}^{(k)}R^{(k+M-1)}R^{(k+M-2)}\dots R^{(k+j+1)}x^j \bm p_N^{(k+j)}(x)
  =x^M L_{\bm\varepsilon}^{(k)}\bm p^{(k)}(x),
\end{equation*}
Substituting $x=x\rme^{-2\pi \rmi \nu/M}$, $\nu=0, 1, \dots, M-1$,
and taking a linear combination, we obtain
\begin{multline}\label{eq:trr-eM-BLP-2}
  H^{(k)}\sum_{\nu=0}^{M-1} w_\nu\bm p^{(k)}(x\rme^{-2\pi \rmi \nu/M})
  +\sum_{j=0}^{M-1}L_{\bm{\varepsilon}^*}^{(k)}R^{(k+M-1)}R^{(k+M-2)}\dots R^{(k+j+1)}x^j\sum_{\nu=0}^{M-1} w_\nu \rme^{-2\pi \rmi \nu j/M} \bm p_N^{(k+j)}(x\rme^{-2\pi \rmi \nu/M})\\
  =x^ML_{\bm\varepsilon}^{(k)}\sum_{\nu=0}^M w_\nu \bm p^{(k)}(x\rme^{-2\pi \rmi \nu/M}),
\end{multline}
where $w_0, w_1, \dots, w_{M-1}$ are some constants.
If there exist some values $x_0, x_1, \dots, x_{N-1}$ and
constants $w_{r, 0}^{(k)}, w_{r, 1}^{(k)}, \dots, w_{r, M-1}^{(k)}$ for $r=0, 1, \dots, N-1$
satisfying
\begin{equation*}
  \sum_{\nu=0}^{M-1} w_{r, \nu}^{(k)}\rme^{-2\pi\rmi \nu j/M}p_N^{(k+j)}(x_r \rme^{-2\pi\rmi \nu/M})=0,\quad
  j=0, 1, \dots, M-1,
\end{equation*}
then \eqref{eq:trr-eM-BLP-2} yields
\begin{equation*}
  H^{(k)}\sum_{\nu=0}^{M-1} w_{r, \nu}^{(k)}\bm p^{(k)}(x_r \rme^{-2\pi\rmi\nu/M})
  =x_r^M L_{\bm\varepsilon}^{(k)}\sum_{\nu=0}^{M-1} w_{r, \nu}^{(k)}\bm p^{(k)}(x_r \rme^{-2\pi\rmi\nu/M}),
\end{equation*}
i.e., $x_r^M$ and $\sum_{\nu=0}^{M-1} w_{r, \nu}^{(k)}\bm p^{(k)}(x_r \rme^{-2\pi\rmi\nu/M})$
are an eigenvalue and an eigenvector of $(H^{(k)}, L_{\bm\epsilon}^{(k)})$.

Let us introduce the variables
\begin{gather*}
  d_0^{(k+M)}\coloneq q_0^{(k)}+\epsilon_0 e_0^{(k)},\quad
  d_n^{(k+M)}\coloneq q_n^{(k+M)}-(1-\epsilon_n) e_{n}^{(k)},\quad n=1, \dots, N-1,\\
  f_n^{(k)}\coloneq q_n^{(k)}+\epsilon_n e_n^{(k)},\quad n=0, 1, \dots, N-2,\quad
  f_{N-1}^{(k)}\coloneq q_{N-1}^{(k)}.
\end{gather*}
Then, in the same way as in the previous section,
we obtain
\begin{gather}
  f_n^{(k)}=q_n^{(k)}+\epsilon_n e_n^{(k)},\quad n=0, 1, \dots, N-2,\quad
  f_{N-1}^{(k)}=q_{N-1}^{(k)},\label{eq:d-heToda-f}\\
  d_0^{(k+M)}=f_0^{(k)},\quad
  d_n^{(k+M)}=
  \begin{dcases*}
    d_{n-1}^{(k+M)}\frac{f_n^{(k)}}{q_{n-1}^{(k+M)}} & if $\epsilon_{n-1}=0$,\\
    q_{n-1}^{(k)}\frac{f_n^{(k)}}{f_{n-1}^{(k)}} & if $\epsilon_{n-1}=1$,
  \end{dcases*}\quad
  n=1, 2, \dots, N-1,\\
  q_n^{(k+M)}=d_n^{(k+M)}+(1-\epsilon_n)e_n^{(k)},\quad n=0, 1, \dots, N-1,\\
  e_n^{(k+1)}=e_n^{(k)}\frac{f_{n+1}^{(k)}}{q_n^{(k+M)}+\epsilon_n e_{n-1}^{(k+1)}},\quad
  n=0, 1, \dots, N-2.\label{eq:d-heToda-e}
\end{gather}
We call this system the \emph{discrete hungry elementary Toda orbits}~\cite{kobayashi2022geb}.

\subsection{Isospectral transformation to an upper Hessenberg matrix}
In the same manner as in the previous section,
we can derive the relations
\begin{gather}
  xp_n^{(k+1)}(x)=p_{n+1}^{(k)}(x)+q_n^{(k)}p_n^{(k)}(x),\label{eq:CT-eM-biorth}\\
  xp_n^{(k+1)}(x)=p_{n+1}^{(k+(\epsilon_n-1)M)}(x)+d_n^{(k)}p_n^{(k)}(x),\\
  xp_n^{(k+1)}(x)=p_{n+1}^{(k+\epsilon_n M)}(x)+f_n^{(k)}p_n^{(k)}(x),\label{eq:CT-M1-biorth}\\
  p_n^{(k)}(x)=p_n^{(k+M)}(x)+e_{n-1}^{(k)}p_{n-1}^{(k+(1-\epsilon_{n-1})M)}(x).\label{eq:GT-M1-biorth}
\end{gather}
Let us introduce the variables
\begin{equation*}
  \hat q_n^{(k)}\coloneq f_n^{(k+\eta_n M)},\quad
  \hat e_n^{(k)}\coloneq e_n^{(k+\eta_{n+1} M)},
\end{equation*}
polynomials
\begin{equation*}
  \hat p_n^{(k)}(x)\coloneq p_n^{(k+\eta_n M)}(x),
\end{equation*}
and the bidiagonal matrices and the vectors of polynomials
same as \eqref{eq:bidiagonal-OPS} and \eqref{eq:vector-OPS}.
Then, the relations \eqref{eq:CT-M1-biorth} and \eqref{eq:GT-M1-biorth} are rewritten as
\begin{equation*}
  \hat{\bm p}^{(k-M)}(x)=\hat L^{(k-M)}\hat{\bm p}^{(k)}(x),\quad
  x\hat{\bm p}^{(k+1)}(x)=\hat R^{(k)}\hat{\bm p}^{(k)}(x)+\bm p_N^{(k)}(x).
\end{equation*}
Hence, $\{\hat{p}_n^{(k)}(x)\}_{n=0}^N$ satisfies the $(M+2)$-terms recurrence relation
\begin{equation*}
  \hat H^{(k)}\hat{\bm p}^{(k)}(x)+\sum_{j=0}^{M-1}\hat L^{(k)}\hat R^{(k+M-1)}\hat R^{(k+M-2)}\dots \hat R^{(k+j+1)}x^j \bm p_N^{(k+j)}(x)=x^M \hat{\bm p}^{(k)}(x),
\end{equation*}
where
\begin{equation*}
  \hat H^{(k)}\coloneq \hat L^{(k)}\hat R^{(k+M-1)}\hat R^{(k+M-2)}\dots \hat R^{(k)}.
\end{equation*}
It is known that the polynomial sequence satisfying this relation
becomes the $(M, 1)$-biorthogonal polynomials~\cite{maeda2017nuh};
there exists a linear functional $\mathcal L^{(k)}\colon \mathbb C[x]\to \mathbb C$
satisfying the $(M, 1)$-biorthogonality condition
\begin{equation*}
  \mathcal L^{(k)}[x^{mM} \hat p_n^{(k)}(x)]=h_n^{(k)} \delta_{m, n} ,\quad
  h_n^{(k)}\ne 0,\quad
  n=0, 1, \dots, N-1, \quad
  m=0, 1, \dots, n.
\end{equation*}
Let $\mu_{k+m}$ be the moment of this linear functional which is the same definition
as \eqref{eq:def-moment}.
Then, the determinant expression of $\{\hat p_n^{(k)}(x)\}_{n=0}^N$ is given by
\begin{equation*}
  \hat p_n^{(k)}(x)=\frac{1}{\tau_n^{(k)}}
  \begin{vmatrix}
    \mu_k & \mu_{k+1} & \dots & \mu_{k+n-1} & \mu_{k+n}\\
    \mu_{k+M} & \mu_{k+M+1} & \dots & \mu_{k+M+n-1} & \mu_{k+M+n}\\
    \vdots & \vdots & & \vdots & \vdots\\
    \mu_{k+(n-1)M} & \mu_{k+(n-1)M+1} & \dots & \mu_{k+(n-1)M+n-1} & \mu_{k+(n-1)M+n}\\
    1 & x & \dots & x^{n-1} & x^n
  \end{vmatrix},\quad
  n=1, 2, \dots, N,
\end{equation*}
where $\tau_n^{(k)}$ is the block Hankel determinant of the moments
\begin{equation*}
  \tau_0^{(k)}\coloneq 1,\quad
  \tau_n^{(k)}\coloneq |\mu_{k+iM+j}|_{i, j=0}^{n-1},\quad
  n=1, 2, 3, \dots.
\end{equation*}

Next, we consider the determinant expression of $\{p_n^{(k)}(x)\}_{n=0}^N$:
\begin{multline*}
  p_n^{(k)}(x)=\frac{1}{\tau_n^{(k-\eta_n M)}}
  \begin{vmatrix}
    \mu_{k-\eta_n M} & \mu_{k-\eta_n M+1} & \dots & \mu_{k-\eta_n M+n-1} & \mu_{k-\eta_n M+n}\\
    \mu_{k+(-\eta_n+1) M} & \mu_{k-(-\eta_n+1)M+1} & \dots & \mu_{k+(-\eta_n+1)M+n-1} & \mu_{k+(-\eta_n+1)M+n}\\
    \vdots & \vdots & & \vdots & \vdots\\
    \mu_{k+(-\eta_n+n-1)M} & \mu_{k+(-\eta_n+n-1)M+1} & \dots & \mu_{k+(-\eta_n+n-1)M+n-1} & \mu_{k+(-\eta_n+n-1)M+n}\\
    1 & x & \dots & x^{n-1} & x^n
  \end{vmatrix},\\
  n=1, 2, \dots, N.
\end{multline*}
Hence, $\{p_n^{(k)}(x)\}_{n=0}^N$ satisfies the $(\epsilon, M)$-biorthogonality relation
for $n=0, 1, \dots, N-1$:
\begin{equation*}
  \mathcal L^{(k)}[x^{mM} p_n^{(k)}(x)]=0,\quad
  m=-\eta_n , -\eta_n+1, \dots, -\eta_n+n-1,
\end{equation*}
and
\begin{equation*}
  \mathcal L^{(k)}[x^{(-\eta_n+n-\epsilon_n(n+1))M} p_n^{(k)}(x)]=(1-2\epsilon_n)^n \frac{\tau_{n+1}^{(k-\eta_{n+1} M)}}{\tau_n^{(k-\eta_n M)}},
\end{equation*}
or
\begin{equation*}
  \mathcal L^{(k)}[x^{(-\eta_n+n) M} p_n^{(k)}(x)]=\frac{\tau_{n+1}^{(k-\eta_{n} M)}}{\tau_n^{(k-\eta_n M)}},\quad
  \mathcal L^{(k)}[x^{(-\eta_n-1) M} p_n^{(k)}(x)]=(-1)^n\frac{\tau_{n+1}^{(k+(-\eta_{n}-1) M)}}{\tau_n^{(k-\eta_n M)}}.
\end{equation*}
The monic polynomial sequence $\{p_n^{(k)}(x)\}_{n=0}^N$ is called
the \emph{$(\epsilon, M)$-biorthogonal Laurent polynomials}.
Then, the relations~\eqref{eq:CT-eM-biorth}--\eqref{eq:GT-M1-biorth} and
\begin{equation*}
  p_n^{(k)}(0)=(-1)^n\frac{\tau_n^{(k-\eta_n M+1)}}{\tau_n^{(k-\eta_n M)}},\quad
  p_{n+1}^{(k+\epsilon_n M)}=(-1)^{n+1} \frac{\tau_{n+1}^{(k+(-\eta_{n+1}+\epsilon_n)M+1)}}{\tau_{n+1}^{(k+(-\eta_{n+1}+\epsilon_n) M)}}=(-1)^{n+1} \frac{\tau_{n+1}^{(k-\eta_{n}M+1)}}{\tau_{n+1}^{(k-\eta_{n} M)}}
\end{equation*}
give
\begin{gather*}
  q_n^{(k)}
  =-\frac{p_{n+1}^{(k)}(0)}{p_{n}^{(k)}(0)}
  =\frac{\tau_n^{(k-\eta_n M)}\tau_{n+1}^{(k-\eta_{n+1}M+1)}}{\tau_{n}^{(k-\eta_n M+1)}\tau_{n+1}^{(k-\eta_{n+1}M)}},\\
  d_n^{(k)}
  =-\frac{p_{n+1}^{(k+(\epsilon_n-1)M)}(0)}{p_{n}^{(k)}(0)}
  =\frac{\tau_n^{(k-\eta_n M)}\tau_{n+1}^{(k+(-\eta_{n}-1)M+1)}}{\tau_{n}^{(k-\eta_n M+1)}\tau_{n+1}^{(k+(-\eta_{n}-1)M)}},\\
  f_n^{(k)}
  =-\frac{p_{n+1}^{(k+\epsilon_n M)}(0)}{p_{n}^{(k)}(0)}
  =\frac{\tau_n^{(k-\eta_n M)}\tau_{n+1}^{(k-\eta_{n}M+1)}}{\tau_{n}^{(k-\eta_nM+1)}\tau_{n+1}^{(k-\eta_{n}M)}},
\end{gather*}
and, since
\begin{gather*}
  \mathcal L^{(k)}[x^{(-\eta_{n+1}+n+1-\epsilon_n(n+1))M}p_{n+1}^{(k)}(x)]=
  \begin{dcases*}
    \frac{\tau_{n+2}^{(k-\eta_{n+1}M)}}{\tau_{n+1}^{(k-\eta_{n+1}M)}} & if $\epsilon_n=0$,\\
    0 & if $\epsilon_n=1$,
  \end{dcases*}\\
  \mathcal L^{(k+M)}[x^{(-\eta_{n+1}+n-\epsilon_n(n+1))M}p_{n+1}^{(k+M)}(x)]=
  \begin{dcases*}
    0 & if $\epsilon_n=0$,\\
    (-1)^{n+1}\frac{\tau_{n+2}^{(k-\eta_{n+1}M)}}{\tau_{n+1}^{(k+(-\eta_{n+1}+1)M)}} & if $\epsilon_n=1$,
  \end{dcases*}\\
  \mathcal L^{(k+(1-\epsilon_n)M)}[x^{(-\eta_n+n-\epsilon_n(n+1))M}p_n^{(k+(1-\epsilon_n)M)}(x)]=
  \begin{dcases*}
    \frac{\tau_{n+1}^{(k+(-\eta_n+1) M)}}{\tau_n^{(k+(-\eta_n+1) M)}}=\frac{\tau_{n+1}^{(k+(-\eta_{n+1}+1) M)}}{\tau_n^{(k+(-\eta_{n+1}+1) M)}} & if $\epsilon_n=0$,\\
    (-1)^n \frac{\tau_{n+1}^{(k+(-\eta_n-1)M)}}{\tau_{n}^{(k-\eta_n M)}}=(-1)^n \frac{\tau_{n+1}^{(k-\eta_{n+1}M)}}{\tau_{n}^{(k+(-\eta_{n+1}+1) M)}}& if $\epsilon_n=1$
  \end{dcases*}
\end{gather*}
hold,
\begin{equation*}
  e_n^{(k)}
  =\frac{\mathcal L^{(k)}[x^{(-\eta_{n+1}+n+1-\epsilon_n(n+1))M}p_{n+1}^{(k)}(x)]-\mathcal L^{(k+M)}[x^{(-\eta_{n+1}+n-\epsilon_n(n+1))M}p_{n+1}^{(k+M)}(x)]}{\mathcal L^{(k+(1-\epsilon_n)M)}[x^{(-\eta_n+n-\epsilon_n(n+1))M}p_n^{(k+(1-\epsilon_n)M)}(x)]}
  =\frac{\tau_n^{(k+(-\eta_{n+1}+1) M)}\tau_{n+2}^{(k-\eta_{n+1}M)}}{\tau_{n+1}^{(k+(-\eta_{n+1}+1)M)}\tau_{n+1}^{(k-\eta_{n+1}M)}}.
\end{equation*}
Further, we have
\begin{equation*}
  \hat q_n^{(k)}=f_n^{(k+\eta_n M)}=\frac{\tau_n^{(k)}\tau_{n+1}^{(k+1)}}{\tau_{n}^{(k+1)}\tau_{n+1}^{(k)}},\quad
  \hat e_n^{(k)}=e_n^{(k+\eta_{n+1} M)}=\frac{\tau_n^{(k+M)}\tau_{n+2}^{(k)}}{\tau_{n+1}^{(k+M)}\tau_{n+1}^{(k)}}.
\end{equation*}

Summarizing the above, we obtain Algorithm~\ref{alg:Hes-bi-to-Hes} of
the isospectral transformation.

\begin{algorithm}[t]
  \caption{Isospectral transformation from Hessenberg--bidiagonal matrix pencil to Hessenberg matrix}
  \label{alg:Hes-bi-to-Hes}
\begin{algorithmic}
  \REQUIRE $\{q_n^{(0)}\}_{n=0}^{N-1}$, $\{q_n^{(1)}\}_{n=0}^{N-1}$, $\dots$, $\{q_n^{(M-1)}\}_{n=0}^{N-1}$,  $\{e_n^{(0)}\}_{n=0}^{N-2}$ and $\bm\epsilon\in\{0, 1\}^{N-1}$ (or bidiagonal matrices $R^{(0)}$, $R^{(1)}$, $\dots$, $R^{(M-1)}$, $L_{\bm\epsilon^*}^{(0)}$ and $L_{\bm\epsilon}^{(0)}$)
  \STATE $\eta_0 \leftarrow 0$
  \FOR{$n=1$ to $N-1$}
  \STATE $\eta_n \leftarrow \eta_{n-1}+\epsilon_{n-1}$
  \ENDFOR
  \FOR{$k=0$ to $(\eta_{N-1}+1)M-1$}
  \FOR{$n=0$ to $N-1$}
  \STATE\IfThenElse{$n<N-1$}
                   {$f_n^{(k)} \leftarrow q_n^{(k)}+\epsilon_n e_n^{(k)}$}
                   {$f_n^{(k)} \leftarrow q_n^{(k)}$}
  \ENDFOR
  \FOR{$n=0$ to $N-1$}
  \STATE\IfThenElse{$n=0$}
                   {$d_n^{(k+M)}\leftarrow f_n^{(k)}$}
                   {\IfThenElse{$\epsilon_{n-1}=0$}
                               {$d_n^{(k+M)}\leftarrow d_{n-1}^{(k+M)}\frac{f_n^{(k)}}{q_{n-1}^{(k+M)}}$}
                               {$d_n^{(k+M)}\leftarrow q_{n-1}^{(k)}\frac{f_n^{(k)}}{f_{n-1}^{(k)}}$}
                   }
  \STATE $q^{(k+M)}_n \leftarrow d_n^{(k+M)}+(1-\epsilon_n)e_n^{(k)}$
  \STATE\IfThen{$n<N-1$}
               {$e_n^{(k+1)} \leftarrow e_n^{(k)}\frac{f_{n+1}^{(k)}}{q_n^{(k+M)}+\epsilon_n e_{n-1}^{(k+1)}}$}
  \ENDFOR
  \ENDFOR
  \FOR{$n=0$ to $N-1$}
  \FOR{$k=0$ to $M-1$}
  \STATE $\hat q_n^{(k)} \leftarrow f_n^{(k+\eta_n M)}$
  \ENDFOR
  \STATE\IfThen{$n<N-1$}
               {$\hat e_n^{(0)} \leftarrow e_n^{(\eta_{n+1}M)}$}
  \ENDFOR
  \ENSURE $\{\hat q_n^{(0)}\}_{n=0}^{N-1}$, $\{\hat q_n^{(1)}\}_{n=0}^{N-1}$, $\dots$, $\{\hat q_n^{(M-1)}\}_{n=0}^{N-1}$ and $\{\hat e_n^{(0)}\}_{n=0}^{N-2}$ (or Hessenberg matrix $\hat H^{(0)}=\hat L^{(0)}\hat R^{(M-1)}\dots\hat R^{(0)}$)
\end{algorithmic}
\end{algorithm}

\subsection{Numerical example}
Let us consider
\begin{gather*}
  L_{\bm\epsilon^*}^{(0)}
  =
  \begin{pmatrix}
    1 & 0 & 0 & 0 & 0 & 0\\
    0 & 1 & 0 & 0 & 0 & 0\\
    0 & 0 & 1 & 0 & 0 & 0\\
    0 & 0 & 0 & 1 & 0 & 0\\
    0 & 0 & 0 & 10 & 1 & 0\\
    0 & 0 & 0 & 0 & 11 & 1
  \end{pmatrix},\\
  R^{(0)}=
  \begin{pmatrix}
    1 & 1 & 0 & 0 & 0 & 0\\
    0 & 2 & 1 & 0 & 0 & 0\\
    0 & 0 & 3 & 1 & 0 & 0\\
    0 & 0 & 0 & 4 & 1 & 0\\
    0 & 0 & 0 & 0 & 5 & 1\\
    0 & 0 & 0 & 0 & 0 & 6
  \end{pmatrix},\quad
  R^{(1)}=
  \begin{pmatrix}
    2 & 1 & 0 & 0 & 0 & 0\\
    0 & 3 & 1 & 0 & 0 & 0\\
    0 & 0 & 4 & 1 & 0 & 0\\
    0 & 0 & 0 & 5 & 1 & 0\\
    0 & 0 & 0 & 0 & 6 & 1\\
    0 & 0 & 0 & 0 & 0 & 7
  \end{pmatrix},\quad
  R^{(2)}=
  \begin{pmatrix}
    3 & 1 & 0 & 0 & 0 & 0\\
    0 & 4 & 1 & 0 & 0 & 0\\
    0 & 0 & 5 & 1 & 0 & 0\\
    0 & 0 & 0 & 6 & 1 & 0\\
    0 & 0 & 0 & 0 & 7 & 1\\
    0 & 0 & 0 & 0 & 0 & 8
  \end{pmatrix},
\end{gather*}
i.e.,
\begin{equation*}
  L_{\bm\epsilon^*}^{(0)}R^{(2)}R^{(1)}R^{(0)}
  =
  \begin{pmatrix}
    6 & 18 8 & 9 & 1 & 0 & 0\\
    0 & 24 & 36 & 12 & 1 & 0\\
    0 & 0 & 60 & 60 & 15 & 1\\
    0 & 0 & 0 & 120 & 90 & 18\\
    0 & 0 & 0 & 1200 & 1110 & 306\\
    0 & 0 & 0 & 0 & 2310 & 1722
  \end{pmatrix},
\end{equation*}
and
\begin{equation*}
  L_{\bm \epsilon}^{(0)}
  =
  \begin{pmatrix}
    1 & 0 & 0 & 0 & 0 & 0\\
    -7 & 1 & 0 & 0 & 0 & 0\\
    0 & -8 & 1 & 0 & 0 & 0\\
    0 & 0 & -9 & 1 & 0 & 0\\
    0 & 0 & 0 & 0 & 1 & 0\\
    0 & 0 & 0 & 0 & 0 & 1
  \end{pmatrix}.
\end{equation*}
Then, Algorithm~\ref{alg:Hes-bi-to-Hes} yields
\allowdisplaybreaks
\begin{gather*}
  \hat R^{(0)}
  =
  \begin{pmatrix}
    8 & 1 & 0 & 0 & 0 & 0\\
    0 & \frac{1045}{196} & 1 & 0 & 0 & 0\\
    0 & 0 & \frac{11783226}{1951015} & 1 & 0 & 0\\
    0 & 0 & 0 & \frac{11202591839}{1537751072} & 1 & 0\\
    0 & 0 & 0 & 0 & \frac{1793288934976}{673133562011} & 1\\
    0 & 0 & 0 & 0 & 0 & \frac{3365490}{23369591}
  \end{pmatrix},\\
  \hat R^{(1)}
  =
  \begin{pmatrix}
    \frac{43}{4} & 1 & 0 & 0 & 0 & 0\\
    0 & \frac{249816}{44935} & 1 & 0 & 0 & 0\\
    0 & 0 & \frac{4459329545}{417182311} & 1 & 0 & 0\\
    0 & 0 & 0 & \frac{281563249429787}{25605158734417} & 1 & 0\\
    0 & 0 & 0 & 0 & \frac{61342417293160530}{164176201497170723} & 1\\
    0 & 0 & 0 & 0 & 0 & \frac{654348548}{340773203}
  \end{pmatrix},\\
  \hat R^{(2)}
  =
  \begin{pmatrix}
    \frac{570}{43} & 1 & 0 & 0 & 0 & 0\\
    0 & \frac{5738006}{988855} & 1 & 0 & 0 & 0\\
    0 & 0 & \frac{2131337471900}{284718590021} & 1 & 0 & 0\\
    0 & 0 & 0 & \frac{417593915190317388}{71923747531523615} & 1 & 0\\
    0 & 0 & 0 & 0 & \frac{5065558609120017904}{2778977782301483047} & 1\\
    0 & 0 & 0 & 0 & 0 & \frac{340773203}{103007824}
  \end{pmatrix},\\
  \hat L^{(0)}
  =
  \begin{pmatrix}
    1 & 0 & 0 & 0 & 0 & 0\\
    \frac{686}{95} & 1 & 0 & 0 & 0 & 0\\
    0 & \frac{17736500}{3269329} & 1 & 0 & 0 & 0\\
    0 & 0 & \frac{92158247808}{19114261985} & 1 & 0 & 0\\
    0 & 0 & 0 & \frac{393943905477395}{312887922561632} & 1 & 0\\
    0 & 0 & 0 & 0 & \frac{448520531195}{11555726719792} & 1
  \end{pmatrix}
\end{gather*}
\allowdisplaybreaks[0]
and
\begin{equation*}
  \hat H^{(0)}
  =\hat L^{(0)}\hat R^{(2)}\hat R^{(1)}\hat R^{(0)}
  =
  \begin{pmatrix}
    1140 & \frac{11898}{49} & \frac{46404}{1867} & 1 & 0 & 0\\[0.5em]
    8232 & \frac{94344}{49} & \frac{581274}{1867} & \frac{1189329}{38368} & 1 & 0\\[0.5em]
    0 & \frac{2240400}{2401} & \frac{109654800}{91483} & \frac{646077099}{1880032} & \frac{146061709}{5496967} & 1\\[0.5em]
    0 & 0 & \frac{8121738240}{3485689} & \frac{26863943637}{17908264} & \frac{2245552524}{12320333} & \frac{20808}{1867}\\[0.5em]
    0 & 0 & 0 & \frac{215519585169}{368025856} & \frac{110756457399}{1076059336} & \frac{44037}{4796}\\[0.5em]
    0 & 0 & 0 & 0 & \frac{886300800}{12585025489} & \frac{107940}{112183}
  \end{pmatrix}.
\end{equation*}
The eigenvalues of both $(L_{\bm\epsilon^*}^{(0)}R^{(2)}R^{(1)}R^{(0)}, L_{\bm \epsilon}^{(0)})$
and $\hat H^{(0)}$ are 3188.27018, 2167.22313, 485.176601, 25.4359335,
1.14466512 and 0.749495172,
which are verified by \texttt{numpy.linalg.eig()}.

\section{Conclusion}
In this paper, we have constructed an isospectral transformation from
a Hessenberg--bidiagonal matrix pencil to a Hessenberg matrix.
Since the obtained algorithms do not contain any subtractions,
if all the input values are positive, then all computations can be performed
without loss of digits for floating point arithmetic.

As a future work, we are interested in an isospectral transformation from
a tridiagonal--tridiagonal matrix pencil to some type of sparse matrix.
The GEVP of tridiagonal--tridiagonal matrix pencils has attracted much attention~\cite{li1994ags,vandebril2009qsa}.
It is known that this GEVP defines so-called \Rii polynomials~\cite{ismail1995goc}
and its associated discrete integrable systems are \Rii chain~\cite{maeda2016gea,spiridonov2000stc} and FST chain~\cite{spiridonov2007idt}.
We will be able to derive the isospectral transformation by using the relation between these
two discrete integrable systems.

\section*{Declaration of competing interest}
There is no competing interest.

\section*{Data availability}
No data was used for research described in the article.

\section*{Acknowledgement}
This work was supported by JSPS KAKENHI Grant Numbers JP19J23445, JP21K13837 and JP19H01792.
This work was partially supported by the joint project
``Advanced Mathematical Science for Mobility Society''
of Kyoto University and Toyota Motor Corporation.

\end{document}